\documentclass[12pt]{article}
\usepackage{fullpage}
\usepackage{graphicx}
\usepackage{float}
\restylefloat{figure}

\newcommand{\qed}{\hbox to 0pt{}\hfill$\rlap{$\sqcap$}\sqcup$\vspace{4mm}}
\def\RR{\hbox{I\kern-.2em\hbox{R}}}
\title{On Nonoscillation of Mixed Advanced-Delay Differential Equations 
with Positive and Negative Coefficients}
\author{Leonid Berezansky $^{1}$ \\Department of Mathematics,
Ben-Gurion University of the Negev, \\
Beer-Sheva 84105, Israel
\\
Elena Braverman $^{2}$ \\
Department of Mathematics and Statistics, University of Calgary, \\
2500 University Drive N.W., Calgary, AB T2N 1N4, Canada \\
and \\
Sandra Pinelas $^{3}$ \\
Universidade dos A\c{c}ores, Departamento de Matem\'{a}tica, \\
R. M\~{a}e de Deus, 
9500-321 Ponta Delgada, Portugal
}
\date{}

\begin{document}
\maketitle

\footnotetext[1]{Partially supported by Israeli Ministry of Absorption}
\footnotetext[2]{Partially supported by the NSERC Research Grant}

\newcommand{\rf}[1]{(\ref{#1})}
\newcommand{\beq}[1]{ \begin{equation}\label{#1} }
\newcommand{\eeq}{\end{equation} }


\begin{abstract}
For a mixed (advanced--delay) differential equation with variable delays 
and coefficients
$$
\dot{x}(t) \pm  a(t)x(g(t)) \mp b(t)x(h(t)) = 0, ~t\geq t_0
$$
where 
$$
a(t)\geq 0,~b(t)\geq 0,~ g(t)\leq t,~ h(t)\geq t
$$
explicit nonoscillation conditions are obtained.
\end{abstract}

{\bf AMS(MOS) subject classification.} 34K11

{\bf Keywords:} mixed differential equations, delay equations, advanced 
equations, positive solutions

\section{Introduction}

Differential equations with delayed in advanced arguments occur in
many applied problems, see \cite{Frish,Gandolfo,James,Asada,Dub,Kaddar},
especially in mathematical economics. There are natural delays in impact
and processing in the economic models, like the policy lags 
\cite{Asada}. Government's stabilization policy can be destabilizing due to
delays in policy response \cite{Asada} 
(if the delays are large enough). An advanced term may,
for example, reflect the dependency on anticipated capital stock
\cite{Dub,Kaddar}. However, apart from the case of constant coefficients,
there are only few results on delayed-advanced (mixed) differential 
equations. The present paper partially fills up this gap. It 
mainly deals with non-oscillation problems: we 
obtain criteria when a nonoscillatory solution with a prescribed tendency
(either nonincreasing or nondecreasing) exists. 

Recently results on oscillation of delay differential
equations (DDE) have taken shape of a developed theory
presented in monographs \cite{p1}-\cite{Agar}.
Most of oscillation criteria for DDE can be extended to equations 
of advanced type (ADE) 
(see \cite{p1}-\cite{p5} and also the recent papers \cite{MJ}-\cite{LS}).  
However, for mixed differential equations (MDE), 
i.e., equations with delay and advanced arguments,
the theory is much less developed.
 
In this paper we consider a mixed differential equation 
\begin{equation}
\label{1ab}
\dot{x}(t) +\delta_1a(t)x(g(t))+\delta_2b(t)x(h(t))=0, ~t\geq t_0
\end{equation}
with variable coefficients $a(t)\geq 0, ~b(t)\geq 0$, one delayed 
($g(t)\leq t$) and one advanced ($h(t)\geq t$) arguments. 
To the best of our knowledge, oscillation of such equations has not been
studied before, except partial cases of  autonomous equations
\cite{p6}-\cite{p10}, equations of the second or higher order  
\cite{p11}-\cite{p13} and equations with constant delays \cite{p14}. 
In \cite{p15} nonoscillation only  of
equation (\ref{1ab}) and higher order equations was considered, where 
$\delta_1$ and $\delta_2$ have the same sign. 
In \cite{p16} the author considers a 
differential equation with a deviating argument 
without the assumption that it is either a delay or an advanced equation.
Hence the results of \cite{p16} can be applied to MDE (\ref{1ab}).
The results of the present paper and of \cite{p16} are independent.  

We consider equation (\ref{1ab}) under the following conditions: 
\begin{description}
\item{(a1)} 
$a(t), b(t), g(t), h(t)$ 
are Lebesgue measurable locally essentially
bounded functions, \\
$a(t)\geq 0,$ $b(t)\geq 0;$

\item{(a2)} $g(t)\leq t, ~h(t)\geq t,
~\lim_{t\rightarrow\infty}g(t)=\infty$.
\end{description}

For  equation (\ref{1ab}) we can consider the same initial value problem as for 
delay equations:
\begin{equation}
\label{1a}
x(t)=\varphi(t), ~ t<t_0,~ x(t_0)=x_0.
\end{equation}

\noindent
{\bf Definition.} An absolutely continuous on each interval $[t_0,b]$
function $x:\RR \rightarrow \RR$ is called 
{\em a solution of problem} (\ref{1ab})-(\ref{1a}),
if it satisfies equation (\ref{1ab}) for almost all $t\in [t_0,\infty)$ and
equalities (\ref{1a}) for $t\leq t_0$.
\vspace{2mm}

In the present paper we will not discuss existence and uniqueness 
conditions for a solution of the problem (\ref{1ab})-(\ref{1a}), instead 
as mentioned before we will only discuss asymptotic properties of the 
solutions. 
\vspace{2mm}

\noindent
{\bf Definition.} 
Solution $x(t),~t_0\leq t<\infty$, of a differential equation 
or inequality  is called {\em nonoscillatory} if there exists
$T$ such that $x(t)\neq 0$ for any $t\geq T$ and {\em oscillatory} 
otherwise.
\vspace{2mm}

In \cite{p17} for  equation (\ref{1ab}) the cases 
when the delay and the advanced 
term have the same sign ($\delta_1=\delta_2=1$ 
or $\delta_1=\delta_2=-1$) were investigated 
and the following main results were obtained.
\vspace{3mm}

\noindent
{\bf Theorem A} {\em Suppose for the equation
\begin{equation}
\label{A1} 
\dot{x}(t)+a(t)x(g(t))+b(t)x(h(t))=0
\end{equation}
 (a1)-(a2) hold, functions  $a(t), b(t), g(t), h(t)$ 
are equicontinuous on $[0, \infty)$ and 
\begin{equation}
\label{01}
\limsup_{t\rightarrow\infty} [t-g(t)]<\infty,
\quad \limsup_{t\rightarrow\infty} [h(t)-t)]<\infty.
\end{equation}
If the delay equation
$$
 \dot{x}(t)+a(t)x(g(t))+b(t)x(t)=0
$$
has a nonoscillatory solution then equation (\ref{A1}) 
also has a nonoscillatory solution.

In particular, if
$$
\limsup_{t\rightarrow\infty} \int_{g(t)}^t a(s)\exp\left\{
\int_{g(s)}^s b(\tau)d\tau\right\} ds<\frac{1}{e}
$$
then equation (\ref{A1}) has a nonoscillatory solution.
}
\vspace{3mm}

\noindent
{\bf Theorem B} {\em Suppose for the equation
\begin{equation}
\label{A2} 
\dot{x}(t)-a(t)x(g(t))-b(t)x(h(t))=0
\end{equation}
(a1)-(a2) hold, 
functions  $a(t), b(t), g(t), h(t)$ are equicontinuous 
on $[0, \infty)$ and condition (\ref{01})
holds.
If the advanced equation
$$
 \dot{x}(t)-a(t)x(t)-b(t)x(h(t))=0
$$
has a nonoscillatory solution then equation (\ref{A2}) also has a nonoscillatory solution.

In particular, if
$$
\lim_{t\rightarrow\infty}\sup \int_t^{h(t)} b(s)\exp\left\{
\int_s^{h(s)} a(\tau)d\tau\right\} ds<\frac{1}{e}
$$
then equation (\ref{A2}) has a nonoscillatory solution.
}
\vspace{2mm}

In the present paper we consider the two remaining cases where 
coefficients have different signs:
$\delta_1=1,~\delta_2=-1$ and $\delta_1=-1,~\delta_2=1$.
For these cases we obtain some rather natural explicit  
nonoscillation conditions for  equation (\ref{1ab}): in Section 2 
for the former case, in Section 3 for the latter one. Finally, Section 4
involves discussion of the results and outlines some open problems on MDE.

It is important to emphasize that in applications of ADE and MDE, for 
example in economics, it is interesting to obtain not only positive 
solutions but also positive monotone solutions, which keep the trend.
Most of the theorems in this paper present sufficient conditions when 
solutions of this kind  exist. 
Some of the results also present explicit estimates for positive 
solutions.

We note that the methods applied in the present paper are different from 
\cite{p17}. The criteria (like Theorem 2) and non-explicit results of the 
general form (Theorem 3) are supplemented by corollaries which provide
easily verified sufficient conditions.

\section{Positive Delay Term, Negative Advanced Term}

In this section we consider the case $\delta_1=1,~\delta_2=-1$ in 
(\ref{1ab}) which becomes
\begin{equation}
\label{1}
\dot{x}(t)+a(t)x(g(t))- b(t)x(h(t))=0, ~~t\geq t_0.
\end{equation}

\newtheorem{guess}{Theorem}
\begin{guess}
Suppose (a1)-(a2) hold and $a(t)\geq b(t)$.
Then the following conditions are equivalent:
\noindent

1. Differential inequality
\begin{equation}
\label{2}
\dot{x}(t)+a(t)x(g(t))- b(t)x(h(t))\leq 0, ~~t\geq t_0,
\end{equation}
has an eventually nonincreasing positive solution.
\noindent

2. Integral inequality
\begin{equation}
\label{3}
u(t)\geq a(t)
\exp\left\{ \int_{g(t)}^t u(s)ds\right\} 
-b(t)\exp\left\{-\int_t^{h(t)} u(s)ds\right\}, ~~t\geq t_1,  
\end{equation}
has a nonnegative locally integrable solution for some $t_1\geq t_0$,
where we assume $u(t)=0$ for $t<t_1$. 

3. Differential equation (\ref{1}) 
has an eventually positive nonincreasing solution.
\end{guess}
{\bf Proof.} $1)\Rightarrow 2)$. 
Let $x$ be a solution of (\ref{2}) such that $x(t)>0,
~~\dot{x}(t)\leq 0,~~t\geq t_0$.
For some $t_1\geq t_0$ we have $g(t)\geq t_0$ for $t\geq t_1$.
Denote $u(t)=-\dot{x}(t)/x(t),~t\geq t_1$, $u(t)=0$, $t<t_1$.
Then 
\begin{equation}
\label{4}
x(t)=x(t_1)\exp\left\{ -\int_{t_1}^t u(s)ds\right\} ,~~t\geq t_1.
\end{equation}
After substituting (\ref{4}) into (\ref{2}) and carrying the
exponent
out of the brackets we obtain
$$
-\exp\left\{ -\int_{t_1}^t u(s)ds\right\}
x(t_1)\left
[u(t)- a(t)\exp \left\{ \int_{g(t)}^t u(s)ds \right\} + b(t)
\exp\left\{ - \int_t^{h(t)} u(s)ds \right\} \right]\leq 0.
$$
Hence (\ref{3}) holds.

$2)\Rightarrow 3)$.
Suppose $u_0(t)\geq 0,~t\geq t_1$ is a solution of inequality (\ref{3}). 
Consider the following sequence
\begin{equation}
\label{5}
u_{n+1}(t) = a(t) \exp \left\{ \int_{g(t)}^t u_n(s)ds \right\} - 
b(t) \exp\left\{ -\int_t^{h(t)} u_n(s)ds \right\}, n\geq 0.
\end{equation}
Since $u_n(t)\geq a(t)-b(t)\geq 0$ and $u_0\geq u_1$ then by induction  
$$
0\leq u_{n+1}(t) \leq  u_n(t)\leq \dots \leq u_0(t).
$$
Hence there exists a pointwise limit 
$$
u(t)=\lim_{n\rightarrow\infty}u_n(t).
$$
The Lebesgue convergence theorem and (\ref{5}) imply 
$$
u(t)=a(t)\exp\left\{ \int_{g(t)}^t u(s)ds \right\} -
b(t)\exp\left\{ -\int_t^{h(t)} u(s)ds \right\}.
$$
Then $x(t)$ denoted by (\ref{4}) 
is a nonnegative nonincreasing solution of  equation (\ref{1}).

Implication $3)\Rightarrow 1)$ is evident. 
\qed

For comparison consider now the following MDE
\begin{equation}\label{6}
\dot{x}(t)+a_1(t)x(g_1(t))- b_1(t)x(h_1(t))=0, ~~t\geq t_0.
\end{equation}

\noindent
{\bf Corollary 1.1}
Suppose (a1)-(a2) hold for $a$,$b$,$h$,$g$,$a_1$,$b_1$,$h_1$,$g_1$ and  
\begin{equation}
\label{11a}
 b_1(t)\leq b(t) \leq a(t)\leq a_1(t),
~ g(t)\geq g_1(t), ~h(t)\leq h_1(t).
\end{equation}
If  equation (\ref{6}) has an eventually positive solution 
with an eventually nonpositive derivative
then the same is valid for  equation (\ref{1}).
\\
{\bf Proof.} 
Suppose (\ref{6}) has an eventually positive solution with an eventually nonpositive 
derivative. By Theorem 1 the integral inequality
$$
u(t)\geq a_1(t)
\exp\left\{ \int_{g_1(t)}^t u(s)ds\right\} 
-b_1(t)\exp\left\{-\int_t^{h_1(t)} u(s)ds\right\}, ~~t\geq t_1,  
$$
has a nonnegative locally integrable solution $u(t)$ for some $t_1$.
Inequalities (\ref{11a}) imply that $u(t)$ also 
satisfies (\ref{3}). Thus by Theorem 1  equation (\ref{1}) has 
an eventually positive solution with an eventually nonpositive
derivative. 
\qed

\noindent
{\bf Corollary 1.2}
Suppose (a1)-(a2) hold,  for $t$  sufficiently large $a(t)\geq b(t)$ and 
$$
b(t)\geq a(t) \left[ \exp \left\{\int_{g(t)}^t a(s)ds \right\}-1\right]
\exp \left\{\int_t^{h(t)} a(s)ds \right\} .
$$
Then there exists an eventually positive solution with an eventually nonpositive
derivative of  equation (\ref{1}).
\\
{\bf Proof}. It is easy to see that $u(t)=a(t)$ is a nonnegative 
solution of inequality (\ref{3}).
\qed

\noindent   
{\bf Corollary 1.3}
Suppose (a1)-(a2) hold and 
there exist $a>0 ,b>0 ,\tau>0 ,\sigma >0$ such that
$$
b\leq b(t)\leq a(t)\leq a,~ g(t)\geq t-\tau, ~ h(t)\leq t+\sigma.
$$
If there exists a solution  $\lambda>0$ of the algebraic equation
\begin{equation}
\label{11}
-\lambda+ae^{\lambda \tau}-be^{-\lambda \sigma}=0
\end{equation}  
then equation (\ref{1}) has an eventually positive solution with an eventually nonpositive
derivative. 
\\
{\bf Proof.} The function $x(t)=e^{-\lambda t}$ is a positive solution of the 
equation
\begin{equation}
\label{10}
\dot{x}(t)+ax(t-\tau)-bx(t+\sigma)=0, t\geq 0.
\end{equation}
By Corollary 1.1  equation (\ref{1}) has a nonoscillatory solution.
\qed

\noindent   
{\bf Corollary 1.4}
Suppose (a1)-(a2) hold, $a(t)\geq b(t)$ and 
there exists a nonoscillatory solution of the delay equation
\begin{equation}
\label{10a}
\dot{x}(t)+a(t)x(g(t))=0.
\end{equation}
Then there exists an eventually positive solution with an eventually nonpositive
derivative  of equation (\ref{1}). 
\\
{\bf Proof.}
Theorem 1 in \cite{p18} implies that there exists a nonnegative solution 
$u(t)$
of the inequality
$$
u(t)\geq a(t)\exp\left\{\int_{g(t)}^t u(s)ds\right\}. 
$$
Hence $u(t)$ is also a nonnegative solution of inequality (\ref{3}),
then by Theorem 1 equation (\ref{1}) has a nonoscillatory solution.
\qed

\noindent
{\bf Remark.} Equation (\ref{10a}) has a nonoscillatory solution
\cite{GL} if for $t$ sufficiently large
$$\int_{g(t)}^t a(s)~ds\leq  \frac{1}{e}.$$

\noindent   
{\bf Corollary 1.5}
Suppose (a1)-(a2) hold, $a(t) \geq b(t)$, the integral inequality
(\ref{3}) has a nonnegative solution for $t \geq t_1$ and 
 $\int_0^{\infty} [a(s)-b(s)]ds=\infty$.
Then there exists 
an eventually positive solution $x(t)$ with an eventually nonpositive
derivative  of equation (\ref{1})
such that $\lim_{t\rightarrow\infty}x(t)=0$.
\\
{\bf Proof.} 
By the assumption of the theorem there exists a nonnegative 
solution $u(t)$ of equation (\ref{3}), which obviously satisfies
$u(t)\geq a(t)-b(t)$. Then the function
$x(t)$ defined by (\ref{4}) is a solution of (\ref{1}).
For this solution we have ${\displaystyle 0<x(t)\leq x(t_1) 
\exp \left\{ -\int_0^t [a(s)-b(s)]ds \right\} }$. 
Hence $\lim_{t\rightarrow\infty}x(t)=0$.
\qed

Corollaries 1.4 and 1.5 imply the following statement.
\vspace{2mm}

\noindent   
{\bf Corollary 1.6}
Suppose (a1)-(a2) hold, $a(t) \geq b(t)$,  for $t$ sufficiently large
$\int_{g(t)}^t a(s)~ds\leq  \frac{1}{e}$ and 
 $\int_0^{\infty} [a(s)-b(s)]\,ds=\infty$.
Then the equation 
$$
\dot{x}(t)+a(t)x(t)-b(t)x(h(t))=0
$$
has an eventually positive solution $x(t)$ with an eventually 
nonpositive derivative  
such that $\lim_{t\rightarrow\infty}x(t)=0$.
\vspace{2mm}

\noindent 
{\bf Example 1.} Consider the equation
\begin{equation}
\label{ex1eq1}
\dot{x}(t) + 1.4 x(t-0.3) - 1.3 x(t+0.3) =0.
\end{equation}
Then $u(t) \equiv 1$ is a solution of the relevant inequality
(\ref{3}) since $1.4e^{0.3}-1.3 e^{-0.3} \approx 0.9267 <1$.
We remark that since $1.3 \cdot 0.3 = 0.39 > \frac{1}{e} \approx 0.368$, 
then all solution of both equations
\begin{equation}
\label{ex1eq2}
\dot{x}(t) + 1.4 x(t-0.3)  =0
\end{equation}
and
\begin{equation}
\label{ex1eq3}
\dot{x}(t)  - 1.3 x(t+0.3) =0
\end{equation}
are oscillatory \cite{LLZ}. 
The characteristic equation 
\begin{equation}
\label{ex1eq4}
-\lambda + 1.4 e^{0.3 \lambda} - 1.3 e^{-0.3 \lambda}=0
\end{equation}
has three real roots: $\lambda_1 \approx -4.2282$, $\lambda_2 \approx 
0.5436$ and $\lambda_3 \approx 3.3541$, thus ${\displaystyle
e^{-\lambda_1 t}}$, ${\displaystyle e^{-\lambda_2 t}}$ and ${\displaystyle 
e^{-\lambda_3 t}}$ are three nonoscillatory solutions of (\ref{ex1eq1}), 
the first one is unbounded on $[0,\infty)$ and the other two are bounded
and have a negative derivative.
\vspace{2mm}

\noindent
{\bf Example 2.}
Consider the equation
\begin{equation}
\label{ex2eq1}
\dot{x}(t) + (1.375+0.025 \sin~t) x(t-0.3) - (1.325+ 0.025 \cos~t) 
x(t+0.3) =0.
\end{equation}
Since $1.3 \leq 1.325+ 0.025 \cos~t \leq 1.35  \leq 1.375+ 0.025 \sin~t
\leq 1.4$ then applying Corollary 1.1 and the results of Example 1 we 
obtain that (\ref{ex2eq1}) has an eventually positive solution $x(t)$ with 
an eventually nonpositive derivative.
Moreover, since the integral of a continuous nonnegative periodic function
$\int_0^{\infty} (0.05 + 0.025 \sin t - 0.025 \cos t)~dt $  diverges, then 
by Corollary 1.5 this solution satisfies $\lim_{t \to \infty} 
x(t)=0$.

\begin{guess}
Suppose (a1)-(a2) hold and $b(t)\geq a(t)$.  
Then the following conditions are equivalent:
\noindent
 
1. Differential inequality
\begin{equation}
\label{7}
\dot{x}(t)+a(t)x(g(t))- b(t)x(h(t))\geq 0, ~~t\geq t_0,
\end{equation}
has an eventually positive solution with an eventually nonnegative 
derivative.
\noindent
 
2. Integral inequality
\begin{equation}
\label{8}
u(t)\geq b(t)\exp\left\{\int_{t}^{h(t)} u(s)ds\right\} -
a(t)\exp\left\{-\int_{g(t)}^t u(s)ds\right\}, ~~t\geq t_1,
\end{equation}
$u(t)=0, t<t_1,$ has a nonnegative locally integrable solution for some $t_1\geq t_0$.
\noindent
 
3. Differential equation (\ref{1}) has an eventually positive 
solution with an eventually nonnegative derivative.
\end{guess}
{\bf Proof.}
$ 1)\Rightarrow 2).$
Let $x$ be a solution of (\ref{7}) such that $x(t)>0,~~\dot{x}(t)\geq 0,~~
t\geq t_0$.
For some $t_1\geq t_0$ we have $g(t)\geq t_0$ for $t\geq t_1$.
Denote $u(t)=\dot{x}(t)/x(t),~t\geq t_1$.
Then
\begin{equation}
\label{8a}
x(t)=x(t_1)\exp\left\{\int_{t_1}^t u(s)ds\right\} ,~~t\geq t_1.
\end{equation}
After substituting (\ref{8a}) into (\ref{7}) and carrying the
exponent
out of the brackets we obtain
$$
\exp\left\{ \int_{t_1}^t u(s)ds\right\}
x(t_1)\left[
u(t)- b(t)\exp\left\{ \int_t^{h(t)} u(s)ds \right\} +
a(t)\exp\left\{ - \int_{g(t)}^t u(s) ds\right\} \right]\geq 0.
$$
Hence (\ref{8}) holds.
 
$2)\Rightarrow 3).$ 
Let $u_0(t)\geq 0$ be a solution of inequality (\ref{8}).
Consider the sequence
\begin{equation}
\label{9}
u_{n+1}(t)= b(t)\exp\left\{ \int_t^{h(t)} u_n(s)ds \right\} -
a(t)\exp\left\{ - \int_{g(t)}^t u_n(s) ds \right\}, n\geq 0.
\end{equation}
Inequalities $u_n(t)\geq b(t)-a(t)$ and $u_0\geq u_1$ imply 
$$
0\leq u_{n+1}(t)\leq  u_n(t)\leq \dots \leq u_0(t).
$$
Hence there exists a pointwise limit
$$
u(t)=\lim_{n\rightarrow\infty}u_n(t).
$$
By the Lebesgue convergence theorem and (\ref{9}) we have
\begin{equation}
\label{9b}
u(t)= b(t)\exp \left\{ \int_t^{h(t)} u(s)ds \right\} -
a(t)\exp\left\{ - \int_{g(t)}^t u(s)ds\right\}.
\end{equation}
Then $x(t)$ defined by (\ref{8a}) 
is a positive solution of  equation (\ref{1}) with 
a nonnegative derivative.
 
Implication $3)\Rightarrow 1)$ is evident.
\qed

The proofs of the following corollaries can be given very similar to those 
given for Theorem 1 and hence we omit most of them.

\noindent
{\bf Corollary 2.1}
Suppose (a1)-(a2) hold for equations (\ref{1}) and (\ref{6}),
$$
a_1(t)\leq a(t)\leq  b(t)\leq b_1(t),
~ g(t)\geq g_1(t), h(t)\leq h_1(t).
$$
If  equation (\ref{6}) 
has an eventually positive solution with an eventually 
nonnegative derivative
then so does  equation (\ref{1}).

\noindent
{\bf Corollary 2.2}
Suppose (a1)-(a2) hold,  for $t$  sufficiently large $b(t)\geq a(t)$ and 
$$
a(t)\geq b(t)\left[ \exp \left\{ \int_t^{h(t)} b(s)ds \right\} -1 
\right] ~\exp \left\{ \int_{g(t)}^t b(s)ds \right\}.
$$
Then there exists  an eventually positive solution 
with an eventually nonnegative  derivative of 
equation (\ref{1}).
\\
{\bf Proof}. It is easy to see that $u(t)=b(t)$ is 
a nonnegative solution of inequality (\ref{8}).
\qed

\noindent
{\bf Corollary 2.3}
Suppose (a1)-(a2) hold and 
there exist $a>0 ,b>0 ,\tau>0 ,\sigma >0$ such that
$$
a\leq a(t)\leq b(t)\leq b,~ g(t)\geq t-\tau,~ h(t)\leq t+\sigma.
$$
If there exists a positive solution  $\lambda>0$ of the algebraic equation
\begin{equation}
\label{11c}
\lambda+ae^{-\lambda \tau}-be^{\lambda \sigma}=0
\end{equation}
then there exists an eventually positive solution 
with an eventually  nonnegative derivative of equation (\ref{1}).

\noindent
{\bf Corollary 2.4}
Suppose (a1)-(a2) hold, $b(t)\geq a(t)$ and 
there exists a nonoscillatory solution of the advanced equation
\begin{equation}
\label{10c}
\dot{x}(t)-b(t)x(h(t))=0.
\end{equation}
Then there exists an eventually positive solution 
with an eventually nonnegative derivative of equation (\ref{1}).

\noindent
{\bf Remark.} If for $t$ sufficiently large
$$\int_t^{h(t)} b(s)~ds\leq  \frac{1}{e},$$
then \cite{LLZ} there exists a nonoscillatory solution of (\ref{10c}). 
\vspace{2mm}

\noindent
{\bf Corollary 2.5}
Suppose (a1)-(a2) hold, $b(t) \geq a(t)$, the integral inequality
(\ref{8}) has a nonnegative solution for $t \geq t_1$ and 
 $\int_0^{\infty} [b(s)-a(s)] \, ds=\infty$.
Then there exists an eventually positive solution 
with an eventually nonnegative  derivative  $x(t)$ of 
(\ref{1}) such that $\lim_{t\rightarrow\infty}x(t)=\infty$.

\section{Negative Delay Term, Positive Advanced Term}

In this section we consider the following  scalar MDE
(\ref{1ab}) with $\delta_1=-1,~ \delta_2=1$:
\begin{equation}
\label{12}
\dot{x}(t) - a(t)x(g(t))+ b(t)x(h(t))=0, ~~t\geq t_0,
\end{equation}

\begin{guess}
\label{the5}
Suppose that $a(t)$, 
$b(t)$ are continuous bounded on $[0,\infty)$ 
functions, $g(t)$ and $h(t)$ are equicontinuous 
on $[0,\infty)$ functions,  there exist positive numbers 
$a_1$, $a_2$, $b_1$, $b_2$, $\tau$, $\sigma$, $t_1$
such that
\begin{equation}
\label{newstar}
a_1\leq a(t)\leq a_2,~ b_1\leq b(t)\leq b_2,~ 
t-g(t)\leq \tau,~ h(t)-t\leq \sigma, t \geq t_1,
\end{equation}
and 
the following algebraic system
\begin{equation}
\label{30}
\left\{\begin{array}{l}a_2e^{y\tau}-b_1e^{-y\sigma}\leq x,\\
b_2e^{x\sigma}-a_1e^{-x\tau}\leq y
\end{array}\right.
\end{equation}
has a solution $x=d_1>0$, $y=d_2>0$.

Then equation (\ref{12}) has a nonoscillatory solution.
\end{guess}
{\bf Proof.} 
 Define the following operator
in the space ${\bf C}[t_1,\infty)$ of 
all bounded continuous on $[t_1,\infty)$ functions
with the usual sup-norm  
$$
(Au)(t)=a(t)\exp\left\{ 
-\int_{g(t)}^tu(s)ds\right\} -b(t)\exp\left\{ \int_t^{h(t)}u(s)ds\right\}, ~t\geq t_1,
$$
$u(t)=0, t\leq t_1$.
Let $x=d_1$, $y=d_2$ be a positive solution of system (\ref{30}).
 Then the inequality
$-d_2\leq u(t)\leq d_1$  implies $-d_2\leq (Au)(t)\leq d_1$.
It means that $A(S) \subset S$, where $$S=\{u(t):-d_2\leq u(t)\leq 
d_1\}.$$

Now we will prove that  $A(S)$ is a compact set in the space
${\bf C}[t_1,\infty)$. 
Denote the integral operators
$$
(Hu)(t):=\int_{g(t)}^t u(s)ds,~~(Ru)(t):=\int_{t}^{h(t)} u(s)ds.
$$
We have for $u\in S$
$$|(Hu)(t)|\leq \max\{d_1,d_2\}\tau,~|(Ru)(t)|\leq \max\{d_1,d_2\}\sigma.$$
Hence the sets $H(S)$ and $R(S)$ are bounded in the space
${\bf C}[t_1,\infty)$.

Let $u\in S$. Then 
$$
|(Hu)(\tau_2)-(Hu)(\tau_1)|\leq 
\left| \int_{g(\tau_1)}^{g(\tau_2)} |u(s)|ds \right|
+\left| \int_{\tau_1}^{\tau_2}|u(s)|ds \right|
$$
$$
\leq \max\{d_1,d_2\}\left( \left|
g(\tau_2)-g(\tau_1)\right|+\left|\tau_2-\tau_1\right|\right)$$
and, similarly,
$$
|(Ru)(\tau_2)-(Ru)(\tau_1)|\leq \max\{d_1,d_2\}(|h(\tau_2)-h(\tau_1)|+|\tau_2-\tau_1|).
$$
Since $g$, $h$ are equicontinuous in $[t_1,\infty)$, then
functions in $H(S)$ and $R(S)$ are also equicontinuous.
Then the sets  $H(S)$ and $R(S)$ are compact, consequently, 
$A(S)$ is also a compact set.

Schauder's fix point theorem implies that 
there exists a solution $u\in S$ of operator equation
$u=Au$. Therefore the function
${\displaystyle x(t)=x(t_1)\exp\left\{\int_{t_1}^t u(s)ds\right\} }$, 
$t\geq t_1$, is a positive solution of equation (\ref{12}).
\qed

\noindent
{\bf Corollary 3.1}
Suppose that $a(t)$, 
$b(t)$ are continuous bounded on $[0,\infty)$ functions, 
$g(t)$ and $h(t)$ are equicontinuous on $[0,\infty)$ and there exist 
positive numbers 
$a_1$, $a_2$, $b_1$, $b_2$, $\tau$, $\sigma$
such that (\ref{newstar}) is satisfied
and at least one of the following conditions holds:
\vspace{2mm}

1) 
${\displaystyle
b_2<a_1, ~0<a_2-b_1<\frac{1}{\tau+\sigma}\ln\frac{a_1}{b_2},
}$ 

\vspace{2mm}

2)
${\displaystyle
a_2<b_1, ~0<b_2-a_1<\frac{1}{\tau+\sigma}\ln\frac{b_1}{a_2}.
}$
\\
\\
Then equation (\ref{12}) has a nonoscillatory solution.
\\
{\bf Proof.} It suffices to prove that system (\ref{30}) has a positive 
solution.
Suppose condition 1) holds. For the second condition the proof is similar.

Let us define the functions ${\displaystyle 
f(x)=b_2e^{x\sigma}-a_1e^{-x\tau}}$ and  ${\displaystyle 
g(y)=a_2e^{y\tau}-b_1e^{-y\sigma} }$ which are both monotone increasing.
We have $f(0)<0$, $f(x_1)=0$, where 
$$x_1=\frac{1}{\tau+\sigma}\ln\frac{a_1}{b_2}.$$
Since $g(y)$ is a monotone function 
then there exists the monotone increasing inverse function $g^{-1}(x)$, 
for which we have $g^{-1}(x_2)=0$, where $x_2=a_2-b_1>0$. 

Denote ${\displaystyle h(x)=g^{-1}(x)-f(x)}$. Condition 1) implies 
$x_2<x_1$, then $f(x_2)<0$ and thus $h(x_2)>0$.

From the equality $g(y)=a_2e^{y\tau}-b_1e^{-y\sigma}=x$ 
we have $a_2e^{y\tau}\leq x+b_1$ 
 for $y \geq 0$.
Then 
$$g^{-1}(x)\leq \frac{1}{\tau}\ln\frac{x+b_1}{a_2}
\mbox{~~~and~~~} h(x) \leq \frac{1}{\tau} \ln\frac{x+b_1}{a_2} - b_2 
e^{x\sigma}+a_1$$
for $x$ large enough. 
Hence $\lim_{x\rightarrow\infty}h(x)=-\infty.$

Since  $h$ is continuous, 
then there exists $x_0>x_2>0$ such that $h(x_0)=0$. It means that 
$f(x_0)=g^{-1}(x_0)$.
Therefore $x_0, y_0=f(x_0)$ is a solution of the system  (\ref{30}).
Hence equation (\ref{12}) has a nonoscillatory solution.
\qed

If the conditions of Corollary 3.1 do not hold we can apply 
numerical methods to prove that system  (\ref{30}) has a positive 
solution.
\vspace{2mm}

\noindent
{\bf Example 3.} Consider the equation
\begin{equation}
\label{ex3eq1}
\dot{x}(t) - (1.3+0.1 \sin t)x(t-0.1-0.1 \cos t) +
(1.7+0.1\cos t)x(t+0.2+0.1 \sin t)=0, ~~t\geq 0.
\end{equation}
Then $a_1=1.2$, $a_2=1.4$, $b_1=1.6$, $b_2=1.8$, $\tau=0.2$, $\sigma=0.3$
and
(\ref{30}) has a positive solution $x=2$, $y=3$, since
$$a_2e^{y\tau}-b_1e^{-y\sigma}=1.4e^{0.6}-1.6e^{-0.9} \approx 1.9 <
x=2,
$$ $$
b_2e^{x\sigma}-a_1e^{-x\tau}\leq y=1.8e^{0.6}-1.2e^{-0.4} \approx 2.48 < 
y=3.
$$
Hence equation (\ref{ex3eq1}) has a nonoscillatory solution.

\begin{figure}[ht]
\centering
\includegraphics[scale=0.45]{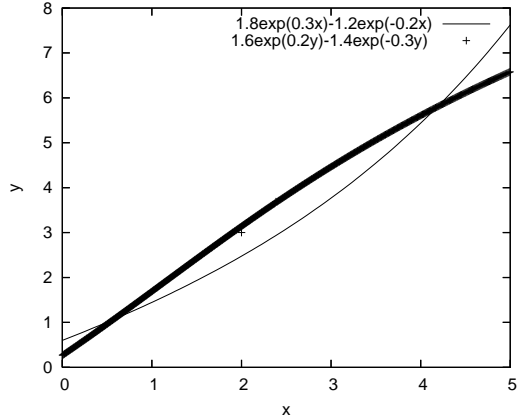}
\caption{The domain of values $(x,y)$ satisfying the system of inequalities
(\protect{\ref{30}}) for equation (\protect{\ref{ex3eq1}}) is between
the curves. The chosen value $x=2$, $y=3$ inside the domain is also marked 
on the graph.
}
\label{figure1}
\end{figure}

Fig. \ref{figure1} illustrates the domain of values $(x,y)$ satisfying the 
system of inequalities (\ref{30}) for equation (\ref{ex3eq1}) which 
is between the two curves.

\vspace{2mm}

\begin{figure}[ht]
\centering
\includegraphics[scale=0.45]{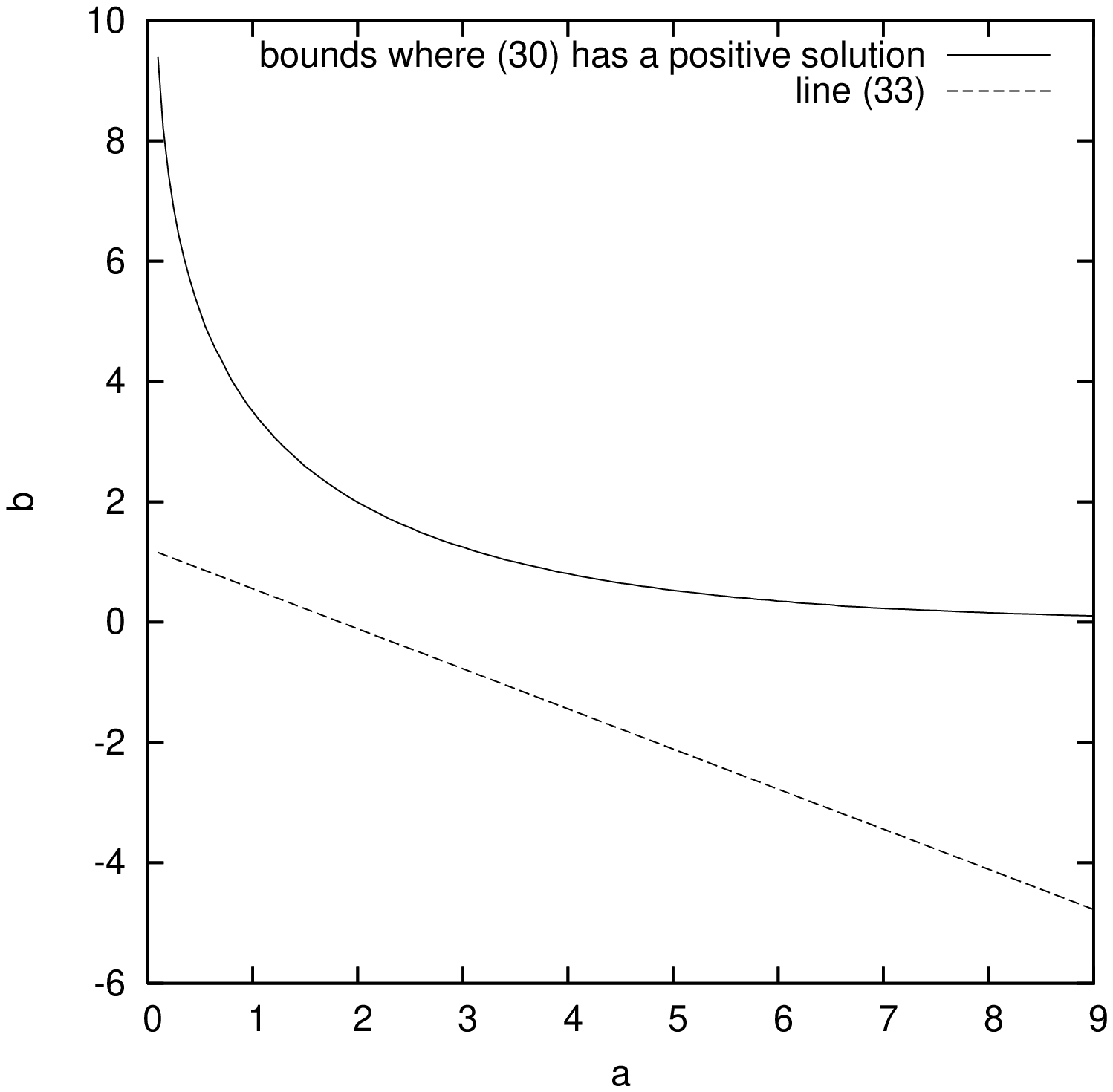}
\caption{The domain of values $a,b$ such that the system of inequalities
(\protect{\ref{30}}) for equation (\protect{\ref{ex4eq1}})
has a positive solution is under
the curve.
}
\label{figure2}
\end{figure}

\noindent
{\bf Example 4.} Consider the equation with constant coefficients and 
variable advance and delay
\begin{equation}
\label{ex4eq1}
\dot{x}(t) - ax(g(t)) + b x(h(t))=0, ~~t\geq 0,
\end{equation} 
where $t \geq g(t) \geq t-0.2$, $t \leq h(t) \leq t+0.3$, like in Example 
3. Thus
in (\ref{30}) we have $\tau=0.2$, $\delta=0.3$. 
All values below the curve in Fig. \ref{figure2}
are such that the system of inequalities
(\ref{30}) has a positive solution and hence equation (\ref{ex4eq1}) has a nonoscillatory solution.

For comparison, we also included the
line
\begin{equation}
\label{ex4eq2}   
0.2 a+0.3b< \frac{1}{e}.
\end{equation}
\vspace{2mm}

\vspace{2mm}

\noindent
{\bf Remark.} The autonomous equation
\begin{equation}
\label{remauton}
\dot{x}(t)-ax(t-\tau)+bx(t+\sigma)=0, ~a>0, ~b>0,~\tau>0, ~\sigma>0, 
\end{equation}
always has a positive solution $e^{\lambda t}$, where $\lambda$ is a 
solution of the characteristic equation
\begin{equation}
\label{remauton1}
f(\lambda)=\lambda-ae^{-\tau\lambda}+be^{\sigma\lambda}=0.
\end{equation}
Since $\lim_{\lambda \to \pm \infty} f(\lambda)= \pm \infty$,
then there is always a real $\lambda$ satisfying
(\ref{remauton1}). There is a positive solution satisfying $\lim_{t 
\to \infty} x(t)=\infty$ if $b<a$ and a positive solution satisfying 
$\lim_{t \to \infty} x(t)=0$ if $b>a$.

\section{Discussion and Open Problems}

In this paper we presented results for equations with variable arguments 
and coefficients, one delay and one advanced term in the case when 
coefficients have any of four possible sign combinations; 
the results for coefficients of different signs are new.
If the delayed term is positive, we not only claim the existence of
a positive solution but present sufficient conditions under which its
asymptotics can be deduced (i.e., a nonincreasing positive solutions which
tends to zero or a nondecreasing solution which tends to infinity).

Below we present some open problems and topics for research and 
discussion.

\begin{enumerate}
\item
If we consider the characteristic equation of the autonomous equation 
(\ref{10}) it is easy to prove that  equation (\ref{10}) 
has a positive solution for any positive coefficients $a>0, b>0$.

Prove or disprove: \\
If conditions (a1)-(a2) hold, 
$$
a_1\leq a(t)\leq a_2,~ b_1\leq b(t)\leq b_2,~ 
t-g(t)\leq \tau~, h(t)-t\leq \sigma,
$$
then equation  (\ref{1}) has an eventually positive solution.

Consider the same problem  for equation   (\ref{12}).

\item
Prove or disprove: \\
If $a(t) \geq b(t) \geq 0$ and (\ref{1}) has an 
eventually positive solution with an eventually nonpositive derivative 
which tends to zero then all solutions of the ordinary differential 
equation
\begin{equation}
\label{ode}
\dot{x}(t)+[a(t)-b(t)] x(t)=0
\end{equation}
tend to zero as $t \to \infty$. If $b(t) \geq a(t) \geq 0$ and (\ref{1}) 
has an
eventually positive solution with an eventually nonnegative derivative
which tends to infinity as $t \to \infty$ then all solutions of 
(\ref{ode}) tend to 
infinity.
\item
For equation (\ref{1}) obtain sufficient conditions under which there 
exists a positive solution which tends to some $d \neq 0$. Is convergence
of ${\displaystyle \int_0^{\infty} a(s)~ds}$ and ${\displaystyle 
\int_0^{\infty} b(s)~ds}$ sufficient, together with the other conditions 
of Theorem 1?  Consider the same problem for (\ref{12}). We remark that
for first order delay neutral equations with positive and negative 
coefficients in \cite{rath1,rath2} sufficient conditions were found 
under which all solutions either oscillate or tend to zero, and for the
second order neutral delay equation such conditions can be found in the 
recent paper \cite{karpuz}. It would be interesting to 
find conditions which guarantee the existence of a positive solution
(which tends to zero) for these equations. For some sufficient conditions 
for linear models without a neutral part see \cite{BDB}.  
\item
Find sufficient conditions when equation (\ref{12}) 
has a positive nonincreasing solution $x(t)$
such that 
$\lim_{t\rightarrow\infty} x(t)=0$
or find sufficient conditions when equation (\ref{12}) 
has a positive solution $x(t)$
with a nonnegative derivative such that $\lim_{t\rightarrow\infty} 
x(t)=\infty$.
\item
Obtain sufficient conditions when the equation with one term which can
be both advanced and delayed and an oscillating coefficient
\begin{equation}
\label{open2}
\dot{x}(t)+a(t)x(g(t))=0
\end{equation}
has a positive solution. For instance, if $a(t)[t-g(t)]$ is either 
positive or negative for any $t$, then (\ref{open2}) 
can be rewritten in the form (\ref{1}) 
or (\ref{12}), thus some conditions can be deduced from the results of 
the present paper. 
\item
Prove or disprove:\\
Suppose ${\displaystyle a(t)\geq 0,~ b(t)\geq 0,~ 
\int_0^{\infty} [a(t)+b(t)]dt=\infty}$.

If equation (\ref{A1}) has a positive solution 
then this equation is asymptotically stable.

If equation (\ref{A2}) has a positive solution 
then the absolute value of any  nontrivial 
solution  tends to infinity.

\item
Instead of initial condition (\ref{1a}) for MDE we can formulate
the ``final" conditions
\begin{equation}
\label{open3}
x(t)=\varphi(t), ~t >t_0, ~x(t_0)=x_0,
\end{equation}
studying nonoscillation for $t<T\leq t_0$ and asymptotics for $t \to 
-\infty$. How would the results of the present paper change if obtained 
for (\ref{1}),(\ref{open3}) or (\ref{12}),(\ref{open3})?
\end{enumerate}

\end{document}